# BOOSTING FOR HIGH-DIMENSIONAL LINEAR MODELS

By Peter Bühlmann

*ETH Zürich*

We prove that boosting with the squared error loss, $L_2$Boosting, is consistent for very high-dimensional linear models, where the number of predictor variables is allowed to grow essentially as fast as $O(\exp(\text{sample size}))$, assuming that the true underlying regression function is sparse in terms of the $\ell_1$-norm of the regression coefficients. In the language of signal processing, this means consistency for de-noising using a strongly overcomplete dictionary if the underlying signal is sparse in terms of the $\ell_1$-norm. We also propose here an *AIC*-based method for tuning, namely for choosing the number of boosting iterations. This makes $L_2$Boosting computationally attractive since it is not required to run the algorithm multiple times for cross-validation as commonly used so far. We demonstrate $L_2$Boosting for simulated data, in particular where the predictor dimension is large in comparison to sample size, and for a difficult tumor-classification problem with gene expression microarray data.

**1. Introduction.** Freund and Schapire's [11] AdaBoost algorithm for classification has attracted much attention in the machine learning community (cf. [20] and the references therein) as well as in related areas in statistics [1, 13], mainly because of its good empirical performance with a variety of datasets. Boosting methods were originally introduced as multiple prediction schemes, averaging estimated predictions from reweighted data. Later, Breiman [1, 2] noted that the AdaBoost algorithm can be viewed as a gradient descent optimization technique in function space. This important insight opened a new perspective, namely to use boosting methods in contexts other than classification. For example, Friedman [12] developed boosting methods for regression which are implemented as an optimization using the squared error loss function: this is what we call $L_2$Boosting. It is essentially the same as Mallat and Zhang's [19] matching pursuit algorithm in signal processing.

Recently, Efron, Hastie, Johnstone and Tibshirani [10] made a connection for linear models between forward stagewise linear regression (FSLR), which seems closely related to $L_2$Boosting, and the $\ell_1$-penalized Lasso [22] or basis pursuit [5]. Roughly speaking: under some restrictive assumptions on the design matrix of a linear model, FSLR approximately yields the set of all Lasso solutions (when varying over the penalty parameter). This intriguing insight may be useful to get a rough picture about $L_2$Boosting via its relatedness to FSLR: it does variable selection and shrinkage, similar to the Lasso. However, it should be stated clearly that the methods are not the same; an example showing a distinct difference between $L_2$Boosting and the Lasso is presented in Section 4.3. Moreover, we point out in Section 2.1 that FSLR and $L_2$Boosting are different algorithms as well.

As the main result, we prove here that $L_2$Boosting for linear models yields consistent estimates in the very high-dimensional context, where the number of predictor variables is allowed to grow essentially as fast as $O(\exp(\text{sample size}))$, assuming that the true underlying regression function is sparse in terms of the $\ell_1$-norm of the regression coefficients. This result is, to our knowledge, the first about boosting in the presence of (fast) growing dimension of the predictor. Some consistency results for boosting with fixed predictor dimension include [17, 18] as well as [25]. Except for Jiang's [17] result, these authors consider versions of boosting either with $\ell_1$-constraints for the boosting aggregation coefficients or, as in [25], with a relaxed version of boosting which we found very difficult to use in practice due to the nonobvious tuning of the relaxation, that is, how fast the boosting aggregation coefficients should decay. The result by Zhang and Yu [25] may be generalized without too much effort to a setting with increasing dimension of the predictor variable, but their theoretical work includes only a rigorous treatment of the classification problem (besides the above mentioned disadvantage of their relaxed boosting algorithm). We believe that it is mainly for the case of high-dimensional predictors where boosting, among other methods, has a substantial advantage over more classical approaches. Some evidence for this will be given in Section 4.1, and other supporting empirical results have been reported in [3] in the different context of low- or high-dimensional additive models for comparing $L_2$Boosting with more traditional methods such as backfitting or MARS (restricted to additive function estimates). Notably, many real datasets nowadays are of high-dimensional nature. Besides the well-documented good empirical performance of boosting, we identify it here as a method which can consistently recover very high-dimensional, sparse functions.

We may also view our result as a consistency property for de-noising using $L_2$Boosting with a strongly overcomplete dictionary. In contrast to a complete dictionary, for example, Fourier- or wavelet-basis, the strongly overcomplete noisy case is not well understood. Our result yields at least the basic property of consistency.

Besides the theoretical consistency result, we propose here a computationally efficient approach for the tuning parameter in boosting, that is, the number of boosting iterations. We give an easily computable definition of degrees of freedom for $L_2$Boosting, and we then propose its use in the corrected $AIC$ criterion. Unlike cross-validation, our $AIC$-tuning does not require boosting to be run multiple times. This makes the $AIC$-type data-driven boosting computationally attractive: depending on the data, it is sometimes as fast as the very efficient LARS algorithm for the Lasso with tuning by its default ten-fold cross-validation [6, 10].

We demonstrate on some simulated examples how our $L_2$Boosting performs for (low- and) mainly high-dimensional linear models, in comparison to the Lasso, forward variable selection, ridge regression, ordinary least squares and a method which has been designed for high-dimensional regression [14]. We also consider a difficult tumor-classification problem with gene expression microarray data: the predictive accuracy of $L_2$Boosting is compared with four other, commonly used classifiers for microarray data, and we briefly indicate the interpretation of the $L_2$Boosting fit along the lines of a linear model fit.

**2. $L_2$Boosting with componentwise linear least squares.** To explain boosting for linear models, consider a regression model

$$Y_i = \sum_{j=1}^{p} \beta_j X_i^{(j)} + \varepsilon_i, \qquad i = 1, \ldots, n,$$

with $p$ predictor variables (the $j$th component of a $p$-dimensional vector $x$ is denoted by $x^{(j)}$) and a random, mean-zero error term $\varepsilon$. More precise assumptions for the model are given in Section 3.

We first specify a base procedure: given some input data $\{(X_i, U_i); i = 1, \ldots, n\}$, where $U_1, \ldots, U_n$ denote some (pseudo-)response variables which are not necessarily the original $Y_1, \ldots, Y_n$, the base procedure yields an estimated function

$$\hat{g}(\cdot) = \hat{g}_{(\mathbf{X},\mathbf{U})}(\cdot),$$

based on $\mathbf{X} = [X_i^{(j)}]_{i=1,\ldots,n; j=1,\ldots,p}$, $\mathbf{U} = (U_1, \ldots, U_n)^T$. Here, we will exclusively consider the componentwise linear least squares base procedure:

$$\hat{g}_{(\mathbf{X},\mathbf{U})}(x) = \hat{\beta}_{\hat{\mathcal{S}}} x^{(\hat{\mathcal{S}})}, \qquad \hat{\beta}_j = \frac{\sum_{i=1}^{n} U_i X_i^{(j)}}{\sum_{i=1}^{n} (X_i^{(j)})^2}, \qquad j = 1, \ldots, p,$$
(2.1)
$$\hat{\mathcal{S}} = \arg\min_{1 \le j \le p} \sum_{i=1}^{n} (U_i - \hat{\beta}_j X_i^{(j)})^2.$$

Thus, the componentwise linear least squares base procedure performs a linear least squares regression against the one selected predictor variable which reduces the residual sum of squares most.

Boosting using the squared error loss, $L_2$Boosting, has a simple structure. Boosting algorithms using other loss functions are described in [12].

$L_2$BOOSTING ALGORITHM.

*Step* 1 (*initialization*). Given data $\{(X_i, Y_i); i = 1, \ldots, n\}$, apply the base procedure yielding the function estimate

$$\hat{F}^{(1)}(\cdot) = \hat{g}(\cdot),$$

where $\hat{g} = \hat{g}_{(\mathbf{X},\mathbf{Y})}$ is estimated from the original data. Set $m = 1$.

*Step* 2. Compute residuals $U_i = Y_i - \hat{F}^{(m)}(X_i), i = 1, \ldots, n$, and fit the real-valued base procedure to the current residuals. The fit is denoted by $\hat{g}^{(m+1)}(\cdot) = \hat{g}_{(\mathbf{X},\mathbf{U})}(\cdot)$ which is an estimate based on the original predictor variables and the current residuals.

Update

$$\hat{F}^{(m+1)}(\cdot) = \hat{F}^{(m)}(\cdot) + \hat{g}^{(m+1)}(\cdot).$$

*Step* 3 (*iteration*). Increase the iteration index $m$ by one and repeat step 2 until a stopping iteration $M$ is achieved.

$\hat{F}^{(M)}(\cdot)$ is an estimator of the regression function $\mathbb{E}[Y|X = \cdot]$. $L_2$Boosting is nothing other than repeated least squares fitting of residuals (cf. [3, 12]). With $m = 2$ (one boosting step), it has already been proposed by Tukey [23] under the name "twicing." In the nonstochastic context, the $L_2$Boosting algorithm is known as "Matching Pursuit" [19], which is popular in signal processing for fitting overcomplete dictionaries.

It is often better to use small step sizes: we advocate here using the step size $\nu$ in the update of $\hat{F}^{(m+1)}$ in step 2, which then becomes

$$(2.2) \qquad \begin{aligned} \hat{F}^{(1)}(\cdot) &= \nu\hat{g}(\cdot), \\ \hat{F}^{(m+1)}(\cdot) &= \hat{F}^{(m)}(\cdot) + \nu\hat{g}^{(m+1)}(\cdot), \qquad m \geq 1, 0 < \nu \leq 1, \end{aligned}$$

where $\nu$ is constant during boosting iterations and is small, for example, $\nu = 0.1$. The parameter $\nu$ can be seen as a shrinkage parameter or alternatively, describing the step size when updating $\hat{F}^{(m+1)}(\cdot)$ along the function $\hat{g}^{(m+1)}(\cdot)$. Small step sizes (or shrinkage) make the boosting algorithm slower and require a larger number $M$ of iterations. However, the computational slow-down often turns out to be advantageous for better out-of-sample empirical prediction performance; see [3, 12].

2.1. *Forward stagewise linear regression.* $L_2$Boosting with componentwise linear least squares is related to forward stagewise linear regression (FSLR), as pointed out by Efron, Hastie, Johnstone and Tibshirani [10]. FSLR differs from $L_2$Boosting with componentwise linear least squares in the update of the new estimate $\hat{F}_m$: instead of using (2.2) which becomes

$$\hat{F}^{(m+1)}(x) = \hat{F}^{(m)}(x) + \nu \hat{\beta}_{\hat{\mathcal{S}}_{m+1}} x^{(\hat{\mathcal{S}}_{m+1})},$$

where $\hat{\beta}_{\hat{\mathcal{S}}_{m+1}}$ is the least squares estimate when fitting the current residuals against the best predictor variable $x^{(\hat{\mathcal{S}}_{m+1})}$, FSLR updates via

$$\hat{F}^{(m+1)}_{\text{FSLR}}(x) = \hat{F}^{(m)}_{\text{FSLR}}(x) + \nu \operatorname{sign}(\hat{\beta}_{\hat{\mathcal{S}}_{m+1}}) x^{(\hat{\mathcal{S}}_{m+1})}.$$

Note that this description of FSLR is equivalent to the one in [10]. In our limited experience, FSLR has about the same prediction accuracy as $L_2$Boosting with componentwise linear least squares. However, we give here two reasons to favor boosting over FSLR. First, the update in FSLR is not scale-invariant whereas the boosting update is on the scale of the current residuals via the magnitude of the least squares estimate $\hat{\beta}_{\hat{\mathcal{S}}_{m+1}}$. It implies that FSLR is often more sensitive to the choice of $\nu$ than boosting. In particular, in case of an orthogonal linear model, $L_2$Boosting has a uniform approximation property for the soft-threshold estimator over all values of the threshold parameter, whereas this nice property does not hold anymore for FSLR [4]. Second, the number of boosting iterations can be reasonably well estimated via degrees of freedom defined as the trace of a boosting hat-matrix, as to be described in Section 2.2. Defining reasonable degrees of freedom which are simple to compute seems not easily possible for FSLR. This has also been pointed out by Efron, Hastie, Johnstone and Tibshirani ([10], comment after formula (4.11)), and they suggest the computationally intensive bootstrap to cope with this problem.

We emphasize that Efron, Hastie, Johnstone and Tibshirani [10] do not advocate using FSLR in practice. They rather focus on the more interesting LARS algorithm.

2.2. *Stopping the boosting iterations.* Boosting needs to be stopped at a suitable number of iterations, to avoid overfitting. The computationally efficient $AIC_c$ criterion in (2.3) below can be used in our context where the base procedure has linear components.

Our goal here is to assign degrees of freedom for boosting. Denote by

$$\mathcal{H}^{(j)} = \mathbf{X}^{(j)}(\mathbf{X}^{(j)})^T / \|\mathbf{X}^{(j)}\|^2, \qquad j = 1, \ldots, p,$$

the $n \times n$ hat-matrix for the linear least squares fitting operator using the $j$th predictor variable $\mathbf{X}^{(j)} = (X_1^{(j)}, \ldots, X_n^{(j)})^T$ only; $\|x\|^2 = x^T x$ denotes the Euclidean norm for a vector $x \in \mathbb{R}^n$. It is then straightforward to show [3] that the $L_2$Boosting hat-matrix, when using the step size $0 < \nu \leq 1$, equals

$$\mathcal{B}_m = I - (I - \nu \mathcal{H}^{(\hat{\mathcal{S}}_m)})(I - \nu \mathcal{H}^{(\hat{\mathcal{S}}_{m-1})}) \cdots (I - \nu \mathcal{H}^{(\hat{\mathcal{S}}_1)}),$$

where $\hat{\mathcal{S}}_i \in \{1, \ldots, p\}$ denotes the component which is selected in the componentwise least squares base procedure in the $i$th boosting iteration.

Using the trace of $\mathcal{B}_m$ as degrees of freedom, we employ a corrected version of $AIC$ (cf. [16]) to define a stopping rule for boosting:

$$
\begin{aligned}
AIC_c(m) &= \log(\hat{\sigma}^2) + \frac{1 + \text{trace}(\mathcal{B}_m)/n}{1 - (\text{trace}(\mathcal{B}_m) + 2)/n}, \\
\hat{\sigma}^2 &= n^{-1} \sum_{i=1}^{n} (Y_i - (\mathcal{B}_m \mathbf{Y})_i)^2, \qquad \mathbf{Y} = (Y_1, \ldots, Y_n)^T.
\end{aligned}
\tag{2.3}
$$

An estimate for the number of boosting iterations is then

$$\hat{M} = \underset{1 \leq m \leq m_{\text{upp}}}{\arg\min} \ AIC_c(m),$$

where $m_{\text{upp}}$ is a large upper bound for the candidate number of boosting iterations.

**3. Consistency of $L_2$Boosting in high dimensions.** We present here a consistency result for $L_2$Boosting in linear models where the number of predictors is allowed to grow very fast as the sample size $n$ increases. Consider the model

$$
\begin{aligned}
Y_i &= f_n(X_i) + \varepsilon_i, \qquad i = 1, \ldots, n, \\
f_n(x) &= \sum_{j=1}^{p_n} \beta_{j,n} x^{(j)}, \qquad x \in \mathbb{R}^{p_n},
\end{aligned}
\tag{3.1}
$$

where $X_1, \ldots, X_n$ are i.i.d. with $\mathbb{E}|X^{(j)}|^2 \equiv 1$ for all $j = 1, \ldots, p_n$ and $\varepsilon_1, \ldots, \varepsilon_n$ are i.i.d., independent from $\{X_s; 1 \leq s \leq n\}$, with $\mathbb{E}[\varepsilon] = 0$. The case with heteroscedastic $\varepsilon_i$'s and potential dependence between $\varepsilon_i$ and $X_i$ is discussed in Remark 1 below. The number of predictors $p_n$ is allowed to grow with the sample size $n$. Therefore, also the predictor $X_i = X_{i,n}$ and the response $Y_i = Y_{i,n}$ depend on $n$, but we usually ignore this in the notation. The scaling of the predictor variables $\mathbb{E}|X^{(j)}|^2 = 1$ is not necessary for running $L_2$Boosting, but it allows one to identify the magnitude of the coefficients $\beta_{j,n}$ [see also assumption (A1) below].

We make the following assumptions:

(A1) The dimension of the predictor in model (3.1) satisfies

$$p_n = O(\exp(Cn^{1-\xi})), \qquad n \to \infty, \text{ for some } 0 < \xi < 1, 0 < C < \infty.$$

(A2) $\sup_{n \in \mathbb{N}} \sum_{j=1}^{p_n} |\beta_{j,n}| < \infty$.
(A3) $\sup_{1 \leq j \leq p_n, n \in \mathbb{N}} \|X^{(j)}\|_\infty < \infty$, where $\|X\|_\infty = \sup_{\omega \in \Omega} |X(\omega)|$ ($\Omega$ denotes the underlying probability space).
(A4) $\mathbb{E}|\varepsilon|^s < \infty$ for some $s > 4/\xi$ with $\xi$ from (A1).

Assumption (A1) allows for a very large predictor dimension relative to the sample size $n$. Assumption (A2) is an $\ell_1$-norm sparseness condition (it could be generalized to $\sum_{j=1}^{p_n} |\beta_{j,n}| \to \infty$ sufficiently slowly as $n \to \infty$, at the expense of additional restrictions on $p_n$). Even if $p_n$ grows, all predictors may be relevant but most of them contribute only with small magnitudes (small $|\beta_{j,n}|$). Assumption (A2) holds for regressions where the number of effective predictors is finite and fixed: that is, the number of $\beta_{j,n} \neq 0$ is independent of $n$ and finite. Assumption (A3) about the boundedness of the predictor variables can be relaxed at the price of more restrictive growth of $p = p_n$: it suffices that $\sup_{1 \leq j \leq p_n} \mathbb{E}|X^{(j)}|^s < \infty$ for some $s \geq 4$ if $p_n = O(n^\alpha)$ where $\alpha = \alpha(s) > 0$ is a number, depending on the number of existing moments $s$, which converges monotonically to $\infty$ as $s$ increases; that is, any polynomial growth of $p_n$ is allowed if the number of moments $s$ is sufficiently large.

THEOREM 1. *Consider the model* (3.1) *satisfying* (A1)–(A4). *Then, the boosting estimate* $\hat{F}^{(m)}(\cdot) = \hat{F}_n^{(m)}(\cdot)$ *with the componentwise linear base procedure from* (2.1) *satisfies: for some sequence* $(m_n)_{n \in \mathbb{N}}$ *with* $m_n \to \infty$ $(n \to \infty)$ *sufficiently slowly,*

$$\mathbb{E}_X |\hat{F}_n^{(m_n)}(X) - f_n(X)|^2 = o_P(1), \qquad n \to \infty,$$

*where $X$ denotes a new predictor variable, independent of and with the same distribution as the $X$-component of the data* $(X_i, Y_i), i = 1, \ldots, n$.

A proof is given in Section 6. Theorem 1 says that $L_2$Boosting recovers the true sparse regression function even if the number of predictor variables is essentially exponentially increasing with sample size $n$. Notably, no assumptions are needed on the correlation structure of the predictor variables.

For the Lasso, a consistency result for high-dimensional regression has been given by Greenshtein and Ritov [15]. We should keep in mind, though, that the Lasso is a different estimator than $L_2$Boosting, as will be demonstrated empirically with an example in Section 4.3.

REMARK 1. Theorem 1 also holds for possibly heteroscedastic errors $\varepsilon_i$ which are potentially dependent on $X_i$, by assuming $(X_1, Y_1), \ldots, (X_n, Y_n)$ i.i.d. and suitable moment conditions for $Y_i$. For the case with bounded $Y_i$, a proof follows as for Corollary 1 below.

3.1. *Binary classification.* The theory of $L_2$Boosting for binary classification with $Y_i \in \{0,1\}$ can be essentially deduced from squared error regression. Bühlmann and Yu [3] argue why $L_2$Boosting is also a reasonable procedure for binary classification. We can always write

$$\begin{aligned}
Y_i &= f_n(X_i) + \varepsilon_i, \\
f_n(x) &= \mathbb{E}[Y|X=x] = \mathbb{P}[Y=1|X=x], \qquad \varepsilon_i = Y_i - f_n(X_i),
\end{aligned} \qquad (3.2)$$

where the $\varepsilon_1, \ldots, \varepsilon_n$ are independent but heteroscedastic with $\mathbb{E}[\varepsilon_i] = 0$ and $\operatorname{Var}(\varepsilon_i) = f_n(X_i)(1 - f_n(X_i))$. When using $L_2$Boosting, we get an estimate for the conditional probability function $\mathbb{P}[Y=1|X=x]$, and the $L_2$Boosting plug-in classifier (for equal misclassification costs) is given by $\hat{C}_n^{(m)}(x) = \mathbb{I}_{[\hat{F}_n^{(m)}(x) > 1/2]}$.

The proof of Theorem 1 essentially goes through and we get the following:

COROLLARY 1. *Consider a binary classification problem with* $(X_1, Y_1), \ldots, (X_n, Y_n)$ *independent and* $Y_i \in \{0,1\}$ *for all* $i = 1, \ldots, n$. *Assume that* $f_n(x) = \mathbb{P}_n[Y=1|X=x] = \sum_{j=1}^{p_n} \beta_{j,n} x^{(j)}$ *and* (A1)–(A3) *hold. Then, for the* $L_2$*Boosting estimate as in Theorem* 1: *for some sequence* $(m_n)_{n \in \mathbb{N}}$ *with* $m_n \to \infty$ $(n \to \infty)$ *sufficiently slowly,*

$$\mathbb{E}_X |\hat{F}_n^{(m_n)}(X) - f_n(X)|^2 = o_P(1), \qquad n \to \infty,$$
$$\mathbb{P}_{X,Y}[\hat{C}_n^{m_n}(X) \neq Y] - L_{n,\text{Bayes}} = o_P(1), \qquad n \to \infty,$$

*where* $L_{n,\text{Bayes}}$ *denotes the Bayes risk* $\mathbb{E}_X[\min\{f_n(X), 1 - f_n(X)\}]$ *and* $X, Y$ *denote new response and predictor variables, independent of and with the same distribution as the data* $(X_i, Y_i)$, $i = 1, \ldots, n$.

The proof is given in Section 6.

**4. Numerical results.**

4.1. *Low-dimensional regression surface within low- or high-dimensional predictor space.* We consider the model

$$\begin{aligned}
X &\sim \mathcal{N}_p(0, V), \qquad Y = f(X) + \varepsilon, \qquad p \in \{3, 10, 100\}, \\
f(X) &= a(V)(1 + 5X^{(1)} + 2X^{(2)} + X^{(3)}), \qquad \varepsilon \sim \mathcal{N}(0, 2^2),
\end{aligned} \qquad (4.1)$$

where $a(V)$ is a scaling factor. The covariance matrix for the predictor variable $X$ and the factor $a(V)$ are chosen as

$$V = I_p, \qquad a(V) = 1 \qquad (4.2)$$

for uncorrelated predictors; or for block-correlated predictors,

(4.3) $$V = \begin{pmatrix} 1 & b & c & 0 & \ldots & \ldots & \ldots & 0 \\ b & 1 & b & c & 0 & \ldots & \ldots & 0 \\ c & b & 1 & b & c & 0 & \ldots & 0 \\ 0 & c & b & 1 & b & c & \ddots & \vdots \\ \vdots & \ddots & \ddots & \ddots & \ddots & \ddots & \ddots & 0 \\ 0 & \ldots & \ldots & \ldots & 0 & c & b & 1 \end{pmatrix},$$
$$b = 0.677, \quad c = 0.323, \quad a(V) = 0.779.$$

The constant $a(V)$ is such that the signal-to-noise ratio $\mathbb{E}|f(X)|^2/\sigma_\varepsilon^2$ is about the same for both model specifications. The model (4.1) with either specification (4.2) or (4.3) has only three effective predictors plus an intercept, all of them contributing to the regression function with different magnitudes (different coefficients). The sample size is denoted by $n$; that is, we generate $n$ i.i.d. realizations $(X_i, Y_i), i = 1, \ldots, n$, from the model.

We use $L_2$Boosting, using shrinkage factor $\nu = 0.1$ [see (2.2)] and the corrected $AIC$ criterion for stopping the boosting iterations [see (2.3)]. We compare it with the Lasso using ten-fold cross-validation for selecting the penalty parameter (i.e., using the default setting from the lars package in R with ten-fold cross-validation—CRAN [6]), with forward variable selection for optimizing the classical $AIC$ criterion, with ordinary least squares (OLS) without variable selection and with ridge regression using the "oracle" ridge-penalty parameter which minimizes the squared error loss over the simulations; the last cannot be used in practice but serves as an optimistic value for the performance of ridge regression. Table 1 reports in detail the mean squared error $\mathrm{MSE} = \mathbb{E}[(\hat{f}(X) - f(X))^2]$ where $X$ is a new test observation, independent from but with the same distribution as the training data. Figure 1 summarizes one of the settings. All results are based on 50 model simulations.

For the high-dimensional (relative to $n$) settings with $p \in \{10, 100\}$, $L_2$Boosting and the Lasso are clearly best for this model with very few effective predictors (see Table 1). Figure 1 displays the good performance of the corrected $AIC$ criterion in (2.3) for stopping the boosting iterations. A detailed comparison of the "oracle"-stopping rule of $L_2$Boosting which stops at the boosting iteration minimizing the mean squared error (see Table 2) can be made to the results in Table 1. Obviously, the "oracle" rule can only be applied for simulated data. We also include in Table 2 the performance of the Lasso with the "oracle" penalty parameter minimizing the mean squared error. $L_2$Boosting and the Lasso perform similarly when using the "oracle"-tuning parameters (see Table 2), while the differences are somewhat more pronounced when comparing $AIC_c$-stopped $L_2$Boosting with Lasso using ten-fold CV tuning (see Table 1).

We also consider the case when $p$ increases exponentially while $n$ grows only linearly. We focus on the model (4.1) with (4.2). The results are given in Table 3. For both $L_2$Boosting and the Lasso, the mean squared error exhibits only a slow increase as $n$ grows linearly and $p$ grows exponentially; compare also with the results from Table 1 with fixed $n = 20$. For this particular example, the Lasso is better for large $p = 300$ than $L_2$Boosting (but this does not imply a general superiority).

4.2. *High-dimensional regression surface with $\ell_1$-coefficients.* We consider here a regression model which fits into the theory of an adaptive

TABLE 1
*MSE for $L_2$Boosting, Lasso, forward variable selection (fwd.var.sel.), ridge with "oracle" penalty (ridge\*) and OLS in model* (4.1) *with* (4.2) *and* (4.3)

| Method       | (4.2), $p = 3$ | (4.2), $p = 10$ | (4.2), $p = 100$ |
|--------------|----------------|-----------------|------------------|
| $L_2$Boost   | 1.658 (0.192)  | 2.318 (0.238)   | 8.792 (0.640)    |
| Lasso        | 1.290 (0.162)  | 3.112 (0.463)   | 8.080 (0.773)    |
| fwd.var.sel. | 1.499 (0.215)  | 3.648 (0.421)   | 13.551 (1.275)   |
| ridge*       | 1.079 (0.117)  | 4.436 (0.392)   | 25.748 (0.637)   |
| OLS          | 1.103 (0.127)  | 5.674 (0.556)   | —                |
|              | (4.3), $p = 3$ | (4.3), $p = 10$ | (4.3), $p = 100$ |
| $L_2$Boost   | 1.054 (0.104)  | 1.649 (0.181)   | 4.643 (0.239)    |
| Lasso        | 1.163 (0.108)  | 3.007 (0.509)   | 3.453 (0.403)    |
| fwd.var.sel. | 1.206 (0.104)  | 2.893 (0.373)   | 12.685 (0.911)   |
| ridge*       | 0.777 (0.079)  | 2.442 (0.226)   | 20.799 (0.538)   |
| OLS          | 1.103 (0.127)  | 5.674 (0.556)   | —                |

Sample size $n = 20$. Estimated standard errors in parentheses.

TABLE 2
*MSE for $L_2$Boosting ($L_2$Boost\*) and the Lasso (Lasso\*), both with "oracle"-tuning parameter*

| Model            | $L_2$Boost*   | Lasso*        |
|------------------|---------------|---------------|
| (4.2), $p = 3$   | 1.103 (0.127) | 1.103 (0.127) |
| (4.3), $p = 3$   | 0.891 (0.100) | 1.075 (0.117) |
| (4.2), $p = 10$  | 2.193 (0.230) | 2.208 (0.262) |
| (4.3), $p = 10$  | 1.404 (0.114) | 1.378 (0.116) |
| (4.2), $p = 100$ | 7.583 (0.593) | 7.116 (0.603) |
| (4.3), $p = 100$ | 2.995 (0.208) | 2.730 (0.234) |

Model (4.1) with (4.2) and (4.3). Sample size $n = 20$. Estimated standard errors in parentheses.

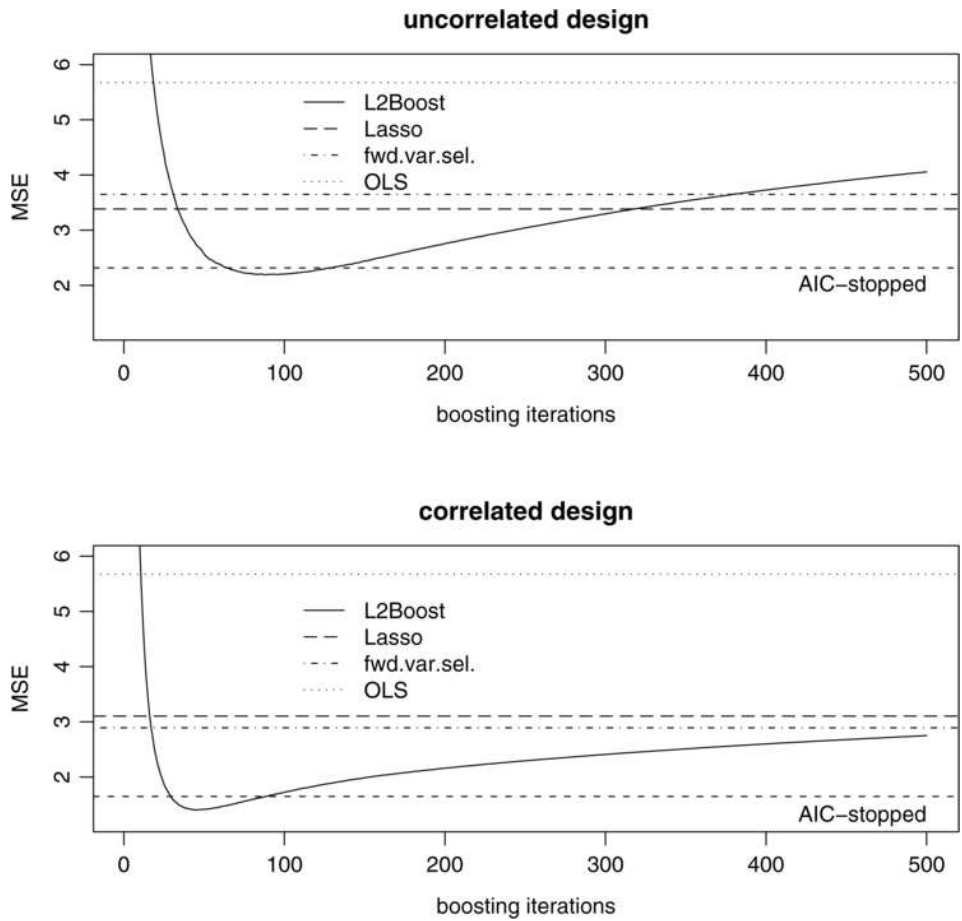

FIG. 1. *MSE for $L_2$ Boosting as a function of boosting iterations (solid line), with $AIC_c$-stopping (dashed line, AIC-stopped), Lasso (long-dashed line), forward variable selection (dashed-dotted line) and OLS (dotted line) in model* (4.1) *with* $p = 10$ *and* (4.2) *(top panel) and* (4.3) *(bottom panel). Sample size* $n = 20$.

estimation procedure for high-dimensional linear regression, presented by Goldenshluger and Tsybakov [14].

The model is

$$
(4.4) \quad \begin{aligned} X &\sim \mathcal{N}_p(0, I), \qquad Y = \sum_{j=1}^{p} \beta_j X^{(j)} + \varepsilon, \\ \beta_j &\sim \mathcal{N}(0, \sigma_j^2), \qquad j = 1, \ldots, p, \qquad \varepsilon \sim \mathcal{N}(0, 1), \end{aligned}
$$

where $\varepsilon, X$ and $\beta_1, \ldots, \beta_p$ are independent of each other. The values $\sigma_j^2$ are decreasing as $j$ increases. Thus, absolute values of the regression coefficients $|\beta_j|$ have a tendency to become small for large $j$. A precise description of the model is given in the Appendix. To summarize, the model is such that $p = p_n$ and $\beta_j = \beta_{j,n}, j = 1, \ldots, p_n$, depend on $n$, satisfying with high probability $\sup_{n \in \mathbb{N}} \sum_{j=1}^{p_n} |\beta_{j,n}| < \infty$, which is our assumption (A2) from Section 3. The sample size is chosen as $n = 100$ and the resulting dimension of the predictor then equals $p = 23$.

We use $L_2$Boosting, using shrinkage $\nu = 0.1$ [see (2.2)] and with estimated number of boosting iterations via the corrected $AIC$ criterion as in (2.3), and we compare it with the Lasso (using the default setting from the lars package in R with ten-fold cross-validation—CRAN [6]), forward variable selection for optimizing the classical $AIC$ criterion, with ridge regression using ten-fold cross-validation for selecting the ridge parameter, with ordinary least squares and with the procedure from [14]. Table 4 displays the results, which are based on 50 independent model simulations. The method from [14] produced one outlier with very large squared error, but the median of the squared errors is still worse than for $L_2$Boosting, Lasso and ridge, which are performing best for this model.

Moreover, the method from [14] depends on the indexing of the predictor variables and is tailored for regression problems where the coefficients $\beta_j$ have a tendency to decay as $j$ increases (e.g., in time series where $j$ indicates the $j$th lagged variable). All other methods do not depend on indexing the predictor variables. We also ran the method from [14] on the same model but with index-reversed regression coefficients

(4.5) $\qquad \beta_1, \ldots, \beta_{23} = \tilde{\beta}_{23}, \ldots, \tilde{\beta}_1, \qquad \tilde{\beta}_j$ as in (4.4).

TABLE 3
*MSE for $L_2$Boosting with $AIC_c$-stopping ($L_2$Boost), with "oracle"-stopping ($L_2$Boost\*), for Lasso with ten-fold CV tuning (Lasso), with "oracle"-tuning (Lasso\*)*

| $(n, p)$  | $L_2$Boost    | Lasso         | $L_2$Boost*   | Lasso*        |
|-----------|---------------|---------------|---------------|---------------|
| (20, 3)   | 1.658 (0.192) | 1.290 (0.162) | 1.103 (0.127) | 1.103 (0.127) |
| (40, 30)  | 2.090 (0.199) | 2.504 (0.274) | 1.730 (0.169) | 1.438 (0.120) |
| (60, 300) | 3.652 (0.186) | 2.136 (0.143) | 2.372 (0.135) | 1.855 (0.122) |

Model (4.1) with (4.2). Estimated standard errors in parentheses.

TABLE 4
*MSE for $L_2$Boosting, Lasso, the method from Goldenshluger and Tsybakov [14] (G&T), forward variable selection (fwd.var.sel.), ridge (ridge) and OLS in model (4.4)*

| $L_2$Boost    | Lasso         | G&T           | fwd.var.sel.  | ridge         | OLS           |
|---------------|---------------|---------------|---------------|---------------|---------------|
| 0.132 (0.006) | 0.135 (0.009) | 0.195 (0.047) | 0.279 (0.019) | 0.116 (0.008) | 0.313 (0.017) |

Sample size $n = 100$. Estimated standard errors in parentheses.

BOOSTING IN HIGH DIMENSIONS 13

The mean squared error (MSE) for the G&T method with (4.5) was then 0.224 (0.025), which shows very clearly the sensitivity of indexing the variables.

4.3. *$L_2$Boosting is different from Lasso.* Consider a model with predictors as in (4.1) and (4.3) with $p = 100$ but with regression function

$$f(X) = 0.2 + 0.2 \sum_{j=1}^{100} X^{(j)} \tag{4.6}$$

and noise $\varepsilon \sim \mathcal{N}(0, 0.5^2)$. The sample size is chosen as $n = 20$. This model is high-dimensional and nonsparse, and it has a high signal-to-noise ratio.

Since all the predictors contribute equally, we may want to keep many of the variables in the model and shrink their corresponding coefficient estimates to zero. However, the Lasso will only allow one to select at most $\min(n, p+1) = 20$ predictor variables (including an intercept); see [26]. When generating one realization of the model (4.6), $L_2$Boosting with the $AIC_c$-stopping rule selected 42 predictor variables (including the intercept), whereas the corresponding number of selected variables with Lasso, tuned by ten-fold cross-validation, is only 13. Thus, we have here an example which demonstrates a feature of $L_2$Boosting which is qualitatively different from the Lasso.

A comparison in terms of the mean squared error yields

$L_2$Boost: 9.468(0.251);    Lasso: 12.140(0.346);    Ridge: 5.548(0.229).

The methods are described in Section 4.2 (estimated standard errors in parentheses). It is no surprise that ridge regression (using ten-fold cross-validation for tuning) performs clearly best. It keeps all variables in the model and shrinks the corresponding estimates toward zero; this is tailored for the structure of the model (4.6) where all the variables contribute equally. We also see from the mean squared error that $L_2$Boosting is quite different (in fact better) than the Lasso. It is not difficult to modify this example such that ridge regression becomes worse than $L_2$Boosting.

4.4. *Gene expression microarray data.* We consider a dataset which monitors $p = 7129$ gene expressions in 49 breast tumor samples using the Affymetrix technology; see [24]. After thresholding to a floor of 100 and a ceiling of 16,000 expression units, we applied a base 10 log-transformation and standardized each experiment to zero mean and unit variance. For each sample, a binary response variable is available, describing the status of lymph node involvement in breast cancer. The data are available at mgm.duke.edu/genome/dna_micro/work/.

TABLE 5
*Cross-validated misclassification rates for lymph node breast cancer data*

|  | $L_2$Boost | FPLR | 1-NN | DLDA | SVM |
|---|---|---|---|---|---|
| Misclassifications | 30.50% | 35.25% | 43.25% | 36.12% | 36.88% |

$L_2$Boosting with *AIC*-stopping ($L_2$Boost), forward variable selection penalized logistic regression (FPLR), 1-nearest-neighbor rule (1-NN), diagonal linear discriminant analysis (DLDA) and a support vector machine (SVM).

We use $L_2$Boosting although the data have the structure of a binary classification problem; Section 3.1 and Corollary 1 yield justification for this, and, for example, Zou and Hastie [26] also use a penalized squared error regression for binary classification with microarray gene expression predictors. The only modification is the *AIC*-stopping criterion: instead of (2.3), we use

$$AIC(m) = -2 \cdot \text{log-likelihood} + 2 \cdot \text{trace}(\mathcal{B}_m),$$

with the Bernoulli log-likelihood. Instead of $L_2$Boosting, we could also use the LogitBoost algorithm [13]: for stopping, the penalty term in the *AIC* criterion above then needs some modification since LogitBoost involves an operator other than $\mathcal{B}_m$.

We estimate the classification performance by a cross-validation scheme where we randomly divide the 49 samples into balanced training- and test-data of sizes $2n/3$ and $n/3$, respectively, and we repeat this 50 times. We compare $L_2$Boosting with *AIC*-stopping (as described above) with four other classification methods: 1-nearest neighbor, diagonal linear discriminant analysis, support vector machine with radial basis kernel (from the R-package `e1071` and using its default values) and a forward selection penalized logistic regression model (using a reasonable penalty parameter and number of selected genes). For 1-nearest neighbor, diagonal linear discriminant analysis and support vector machine, we pre-select the 200 genes which have the best Wilcoxon score in a two-sample problem (estimated from the training dataset only), which is recommended to improve the classification performance; see [9]. Our $L_2$Boosting and the forward variable selection penalized regression are run without pre-selection of genes. The results are given in Table 5. When transforming the response variable to $Y \in \{-1/2, 1/2\}$, that is, subtracting the prior class probability $1/2$, $L_2$Boosting has a cross-validated misclassification rate of 23.13% [4].

For this difficult classification problem, our $L_2$Boosting with componentwise linear least squares (even without centering the response) performs well. It is also interesting to note that the minimal cross-validated misclassification rate as a function of boosting iterations is 29.25%. It shows that the *AIC*-stopping rule is very accurate for this example. A method which we found to perform better for this dataset is the recently proposed Pelora algorithm [7], which does supervised gene grouping: its misclassification rate is 27.88%.

We also show in Figure 2 the estimated regression coefficients for the 42 genes which have been selected during the boosting iterations until *AIC*-stopping; the *AIC* curve is also shown in Figure 2. For comparing the influence of different genes, we display scaled coefficients $\hat{\beta}_j \sqrt{\text{Var}(X^{(j)})}$ which correspond to the estimated coefficients when standardizing the genes to unit variance. There is one gene whose positive expression strongly points toward the class with $Y = 0$ (having negative scaled regression coefficient) and there are five genes whose positive expressions point toward the class with $Y = 1$. The smallest scaled regression coefficient corresponds to a gene which appears as the second best when ranking all the genes with the score of a two-sample Wilcoxon test; the five largest scaled coefficients correspond to the Wilcoxon-based ranks 7, 6, 1, 121, 3 among all the genes. But it should be emphasized that, as usual, our estimated regression model takes partial correlations between the class variable $Y$ and gene expressions (given all other remaining genes) into account, which goes well beyond describing the effects of single genes only.

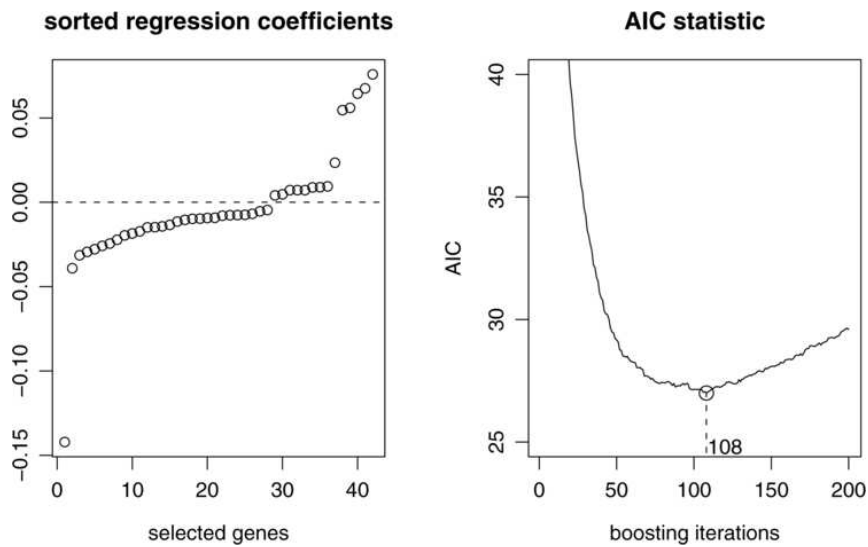

FIG. 2. *Lymph node breast cancer data.* Left: *scaled regression coefficients* $\hat{\beta}_j \sqrt{\text{Var}(X^{(j)})}$ *(in increasing order) from* $L_2$*Boosting for the selected* 42 *genes.* Right: *AIC statistic as a function of* $L_2$*Boosting iterations with minimum at* 108.

**5. Conclusions.** We consider $L_2$Boosting for fitting linear models. The method does variable selection and shrinkage, a property which is very useful in practical applications. This indicates that $L_2$Boosting is related to the $\ell_1$-penalized Lasso, but the methods are not the same.

As a useful device, we propose a simple estimate for the number of boosting iterations, which is the tuning parameter of the method, by using a corrected $AIC_c$ criterion. This makes boosting computationally attractive, since we do not have to run it multiple times in a cross-validation set-up.

We then present some theory for very high-dimensional regression (or for de-noising with strongly overcomplete dictionaries), saying that if the underlying true regression function is sparse in terms of the $\ell_1$-norm of the regression coefficients, $L_2$Boosting consistently estimates the true regression function, even when the number of predictor variables grows like $p_n = O(\exp(n^{1-\xi}))$ for some (small) $\xi > 0$. Notably, no assumptions are made on the correlation structure of the predictors. Thus, we identify $L_2$Boosting as a method which is able, under mild assumptions, to consistently recover very high-dimensional, sparse functions.

**6. Proofs.** We first consider the regression case where the step-size in (2.2) equals $\nu = 1$. In Section 6.3, we give the argument for arbitrary, fixed $0 < \nu \leq 1$. Finally, we present the case for binary classification in Section 6.4.

6.1. *A population version.* The $L_2$Boosting algorithm has a population version which is known as "matching pursuit" [19] or "weak greedy algorithm" [21].

Consider the Hilbert space $L_2(P) = \{f; \|f\|^2 = \int f(x)^2 \, dP(x) < \infty\}$ with inner product $\langle f, g \rangle = \int f(x) g(x) \, dP(x)$. Here, the probability measure $P$ is generating the predictor $X$ in model (3.1). To be precise, the probability measure $P = P_n$ depends on $n$ since the dimensionality of $X$ is growing with $n$: we are actually looking at a sequence of Hilbert spaces $L_2(P_n)$, but we often ignore this notationally [a uniform bound in (6.5) will be a key result to deal with such sequences of Hilbert spaces].

Denote the components of $X$ by
$$g_j(x) = x^{(j)}, \qquad j = 1, \ldots, p_n.$$

Note that by assumption, $\|g_j\| = 1$ for all $j$. Define the following sequence of remainder functions, called matching pursuit or weak greedy algorithm:

(6.1) $$\begin{aligned} R^0 f &= f, \\ R^m f &= R^{m-1} f - \langle R^{m-1} f, g_{\mathcal{S}_m} \rangle g_{\mathcal{S}_m}, \qquad m = 1, 2, \ldots, \end{aligned}$$

where $\mathcal{S}_m$ would be ideally chosen as
$$\mathcal{S}_m = \arg\max_{1 \leq j \leq p_n} |\langle R^{m-1} f, g_j \rangle|.$$

The choice function $\mathcal{S}_m$ is sometimes infeasible to realize in practice. A weaker criterion is: for every $m$ (under consideration), choose any $\mathcal{S}_m$ which satisfies

(6.2) $\quad |\langle R^{m-1}f, g_{\mathcal{S}_m}\rangle| \geq b \cdot \sup_{1 \leq j \leq p_n} |\langle R^{m-1}f, g_j\rangle| \quad$ for some $0 < b \leq 1$.

Of course, the sequence $R^m f = R^{m,\mathcal{S}} f$ depends on $\mathcal{S}_1, \mathcal{S}_2, \ldots, \mathcal{S}_m$, how we actually make the choice in (6.2). Again, we will ignore this notationally.

It easily follows that

$$f = \sum_{j=0}^{m-1} \langle R^j f, g_{\mathcal{S}_{j+1}}\rangle g_{\mathcal{S}_{j+1}} + R^m f$$

and

(6.3) $\quad \|R^m f\|^2 = \|R^{m-1} f\|^2 - |\langle R^{m-1}f, g_{\mathcal{S}_m}\rangle|^2.$

6.1.1. *Temlyakov's result.* Temlyakov [21] gives a uniform bound for the algorithm in (6.1) with (6.2).

If the function $f$ is representable as

(6.4) $\quad f(x) = \sum_j \beta_j g_j(x), \qquad \sum_j |\beta_j| \leq B < \infty,$

which is true by our assumption (A2), then

(6.5) $\quad \|R^m f\| \leq B(1 + mb^2)^{-b/(2(2+b))}, \qquad 0 < b \leq 1$ as in (6.2).

By construction, $R^m f$ depends on the selectors $\mathcal{S}_1, \ldots, \mathcal{S}_m$ in (6.2). The mathematical power of the bound in (6.5) is that it holds for *any* selectors $\mathcal{S}_1, \ldots, \mathcal{S}_m$ which satisfy (6.2). In particular, the bound also holds for sequences $R^m f$ which depend on the sample size $n$ (since $X \sim P = P_n$ and also the function of interest $f = f_n$ depend on $n$).

6.2. *Asymptotic analysis as sample size increases.* The $L_2$ Boosting algorithm can be represented analogously to (6.1). We introduce the notation

$$\langle f, g\rangle_{(n)} = n^{-1} \sum_{i=1}^n f(X_i) g(X_i) \quad \text{and} \quad \|f\|_{(n)}^2 = n^{-1} \sum_{i=1}^n f(X_i)^2$$

for functions $f, g : \mathbb{R}^{p_n} \to \mathbb{R}$. Without loss of generality (but simplifying the notation), we assume in the sequel that $\|g_j\|_{(n)} \equiv 1$ for all $j$ and $n$ (note that $\|g_j\| \equiv 1$ holds already); the justification follows from Lemma 1(i) below. As before, we denote by $\mathbf{Y} = (Y_1, \ldots, Y_n)^T$ the vector of response variables.

Define

$$\hat{R}_n^0 f = f, \qquad \hat{R}_n^1 f = f - \langle \mathbf{Y}, g_{\hat{\mathcal{S}}_1}\rangle_{(n)} g_{\hat{\mathcal{S}}_1},$$
$$\hat{R}_n^m f = \hat{R}_n^{m-1} f - \langle \hat{R}_n^{m-1} f, g_{\hat{\mathcal{S}}_m}\rangle_{(n)} g_{\hat{\mathcal{S}}_m}, \qquad m = 2, 3, \ldots,$$

where
$$\hat{\mathcal{S}}_1 = \argmax_{1 \le j \le p_n} |\langle \mathbf{Y}, g_j \rangle_{(n)}|,$$
$$\hat{\mathcal{S}}_m = \argmax_{1 \le j \le p_n} |\langle \hat{R}_n^{m-1} f, g_j \rangle_{(n)}|, \qquad m = 2, 3, \ldots.$$

By definition, $\hat{R}_n^m f = f - \hat{F}_n^m$ is the difference of the function $f$ and its $L_2$Boosting estimate $\hat{F}_n^m$. Note that we emphasize here the dependence of $\hat{R}_n^m$ on $n$ since finite-sample estimates $\langle \hat{R}_n^{m-1} f, g_j \rangle_{(n)}$ are involved.

6.2.1. *A semipopulation version.* For analyzing $\hat{R}_n^m f$, we want to use Temlyakov's [21] result from (6.5). We will apply it to a semipopulation version $\tilde{R}_n^m f$, as defined below [since it seems difficult to establish (6.2) for $\hat{R}_n^m f$ directly].

Consider
$$\tilde{R}_n^0 f = f,$$
$$\tilde{R}_n^m f = \tilde{R}_n^{m-1} f - \langle \tilde{R}_n^{m-1} f, g_{\hat{\mathcal{S}}_m} \rangle g_{\hat{\mathcal{S}}_m}, \qquad m = 1, 2, \ldots,$$

where $\hat{\mathcal{S}}_m$ is the selector from the sample version above.

The strategy will be as follows. First, we want to establish a finite-sample analogue of (6.2) for the estimated selectors $\hat{\mathcal{S}}_m$; this will then allow us to use Temlyakov's [21] result from (6.5) for $\tilde{R}_n^m f$. Finally, we need to analyze the difference $\hat{R}_n^m f - \tilde{R}_n^m f$.

6.2.2. *Uniform laws of large numbers.*

LEMMA 1. *Under the assumptions* (A1)–(A4), *with* $0 < \xi < 1$ *as in* (A1):

(i) $\sup_{1 \le j, k \le p_n} |n^{-1} \sum_{i=1}^n g_j(X_i) g_k(X_i) - \mathbb{E}[g_j(X) g_k(X)]| = \zeta_{n,1} = O_P(n^{-\xi/2})$,

(ii) $\sup_{1 \le j \le p_n} |n^{-1} \sum_{i=1}^n g_j(X_i) \varepsilon_i| = \zeta_{n,2} = O_P(n^{-\xi/2})$,

(iii) $\sup_{1 \le j \le p_n} |n^{-1} \sum_{i=1}^n f(X_i) g_j(X_i) - \mathbb{E}[f(X) g_j(X)]| = \zeta_{n,3} = O_P(n^{-\xi/2})$,

(iv) $\sup_{1 \le j \le p_n} |n^{-1} \sum_{i=1}^n g_j(X_i) Y_i - \mathbb{E}[g_j(X) Y]| = \zeta_{n,4} = O_P(n^{-\xi/2})$.

PROOF. For assertion (i), denote $M = \sup_j \|g_j(X)\|_\infty$; see assumption (A3). Then Bernstein's inequality yields for every $\gamma > 0$,

$$\mathbb{P}\left[n^{\xi/2} \sup_{1 \le j, k \le p_n} \left| n^{-1} \sum_{i=1}^n g_j(X_i) g_k(X_i) - \mathbb{E}[g_j(X) g_k(X)] \right| > \gamma \right]$$
$$\le p_n^2 2 \exp\left(-\frac{\gamma^2 n^{1-\xi}}{2(\sigma_g^2 + M^2 \gamma n^{-\xi/2})}\right),$$

where $\sigma_g^2$ is an upper bound for $\text{Var}(g_j(X)g_k(X))$ for all $j, k$ (e.g., $\sigma_g^2 = M^4$). Since $p_n^2 = O(\exp(2C(n^{1-\xi})))$, the right-hand side of the inequality above becomes arbitrarily small for $n$ sufficiently large and $\gamma > 0$ large.

For proving assertion (ii), we have to deal with the unboundedness of the $\varepsilon_i$'s in order to apply Bernstein's inequality. Define the truncated variables

$$\varepsilon_i^{\text{tr}} = \begin{cases} \varepsilon_i, & \text{if } |\varepsilon_i| \leq M_n, \\ \text{sign}(\varepsilon_i) M_n, & \text{if } |\varepsilon_i| > M_n. \end{cases}$$

Then for $\gamma > 0$,

$$\mathbb{P}\left[ n^{\xi/2} \sup_{1 \leq j \leq p_n} \left| n^{-1} \sum_{i=1}^n g_j(X_i) \varepsilon_i \right| > \gamma \right]$$

$$\leq \mathbb{P}\left[ n^{\xi/2} \sup_{1 \leq j \leq p_n} \left| n^{-1} \sum_{i=1}^n g_j(X_i) \varepsilon_i^{\text{tr}} - \mathbb{E}[g_j(X) \varepsilon^{\text{tr}}] \right| > \gamma/3 \right]$$

$$+ \mathbb{P}\left[ n^{\xi/2} \sup_{1 \leq j \leq p_n} \left| n^{-1} \sum_{i=1}^n g_j(X_i)(\varepsilon_i - \varepsilon_i^{\text{tr}}) \right| > \gamma/3 \right]$$

$$+ \mathbb{P}\left[ n^{\xi/2} \sup_{1 \leq j \leq p_n} \left| n^{-1} \sum_{i=1}^n \mathbb{E}[g_j(X_i)(\varepsilon_i - \varepsilon_i^{\text{tr}})] \right| > \gamma/3 \right]$$

$$= I + II + III,$$

since $\mathbb{E}[g_j(X)\varepsilon] = \mathbb{E}[g_j(X)]\mathbb{E}[\varepsilon] = 0$, which we use for $III$. We can bound $I$ again by using Bernstein's inequality:

(6.6) $$I \leq p_n 2 \exp\left( -\frac{(\gamma^2/9) n^{1-\xi}}{2(\sigma_g^2 + M_n^2(\gamma/3) n^{-\xi/2})} \right),$$

where $\sigma_g^2$ is an upper bound for $\text{Var}(g_j(X)\varepsilon^{\text{tr}})$ (e.g., $\sup_j \|g_j(X)\|_\infty^2 \mathbb{E}|\varepsilon|^2$). When using

$$M_n = n^{\xi/4},$$

we can make the right-hand side in (6.6) arbitrarily small since $p_n = O(\exp(Cn^{1-\xi}))$; thus, for every $\delta > 0$,

(6.7) $\quad I \leq \delta \quad$ for $n$ sufficiently large, $\gamma$ sufficiently large.

A bound for $II$ can be obtained as follows:

(6.8) $$\begin{aligned} II &\leq \mathbb{P}[\text{some } |\varepsilon_i| > M_n] \leq n \mathbb{P}[|\varepsilon| > M_n] \leq n M_n^{-s} \mathbb{E}|\varepsilon|^s \\ &= O(n^{1-s\xi/4}) = o(1), \quad n \to \infty, \end{aligned}$$

since $s > 4/\xi$ by assumption (A4).

For $III$ we use the bound

(6.9) $$III \leq \mathbb{I}_{[n^{\xi/2} \sup_j |\mathbb{E}[g_j(X)(\varepsilon - \varepsilon^{\text{tr}})]| > \gamma/3]}.$$

Note that by the independence of $\varepsilon$ (and $\varepsilon^{\text{tr}}$) from $g_j(X)$,

$$\mathbb{E}[g_j(X)(\varepsilon - \varepsilon^{\text{tr}})] = \mathbb{E}[g_j(X)]\mathbb{E}[\varepsilon - \varepsilon^{\text{tr}}].$$

Hence, an upper bound is

$$|\mathbb{E}[g_j(X)(\varepsilon - \varepsilon^{\text{tr}})]| \le M|\mathbb{E}[\varepsilon - \varepsilon^{\text{tr}}]|.$$

The latter can be bounded as

$$|\mathbb{E}[\varepsilon - \varepsilon^{\text{tr}}]| \le \left| \int_{|x|>M_n} (\text{sign}(x)M_n - x)\, dP_\varepsilon(x) \right|$$

$$\le \int \mathbb{I}_{[|x|>M_n]}(M_n + |x|)\, dP_\varepsilon(x)$$

$$= M_n \mathbb{P}[|\varepsilon| > M_n] + \int |x|\mathbb{I}_{[|x|>M_n]}\, dP_\varepsilon(x)$$

$$\le M_n^{1-s}\mathbb{E}|\varepsilon|^s + (\mathbb{E}|\varepsilon|^2)^{1/2}(\mathbb{P}[|\varepsilon|>M_n])^{1/2}$$

$$= O(M_n^{1-s}) + O(M_n^{-s/2}) = o(M_n^{-2}) = o(n^{-\xi/2})$$

since $s > 4/\xi > 4$ ($0 < \xi < 1$). Hence, by using (6.9),

$$III = 0 \quad \text{for } n \text{ sufficiently large}, \gamma > 0 \text{ sufficiently large},$$

and together with (6.7) and (6.8), this proves assertion (ii).

Assertion (iii) follows from (i):

$$\sup_{1 \le j \le p_n, n \in \mathbb{N}} \left| n^{-1}\sum_{i=1}^n f(X_i)g_j(X_i) - \mathbb{E}[f(X)g_j(X)] \right|$$

$$\le \sum_{r=1}^{p_n} |\beta_{r,n}| \sup_{1 \le j,k \le p_n} \left| n^{-1}\sum_{i=1}^n g_j(X_i)g_k(X_i) - \mathbb{E}[g_j(X)g_k(X)] \right|$$

$$\le \sum_{r=1}^{p_n} |\beta_{r,n}| \sup_{1 \le j,k \le p_n} \left| n^{-1}\sum_{i=1}^n g_j(X_i)g_k(X_i) - \mathbb{E}[g_j(X)g_k(X)] \right|$$

$$\le \sum_{r=1}^{p_n} |\beta_{r,n}|\zeta_{n,1} = O_P(n^{-\xi/2}).$$

Assertion (iv) follows from (ii) and (iii). □

6.2.3. *Recursive analysis of $L_2$Boosting.* Denote

$$\zeta_n = \max\{\zeta_{n,1}, \zeta_{n,2}, \zeta_{n,3}, \zeta_{n,4}\} = O_P(n^{-\xi/2})$$

which is a bound for all assertions (i)–(iv) in Lemma 1. Also, we denote by $\omega$ a realization of all $n$ datapoints.

..._x

LEMMA 2. *Under the assumptions of Lemma 1, there exists a constant $0 < C_* < \infty$, independent of $n$ and $m$, such that*

$$\sup_{1 \leq j \leq p_n} |\langle \hat{R}_n^m f, g_j \rangle_{(n)} - \langle \tilde{R}_n^m f, g_j \rangle| \leq (5/2)^m \zeta_n C_*$$

*on the set $A_n = \{\omega; |\zeta_n(\omega)| < 1/2\}$.*

*Note that Lemma 1 implies that $\mathbb{P}[A_n] \to 1$, $n \to \infty$. The constant $C_*$ depends on $\sup_{n \in \mathbb{N}} \sum_{j=1}^{p_n} |\beta_{j,n}|$.*

PROOF. We proceed recursively. For $m = 0$, the statement follows directly from Lemma 1(iv). Denote $A_n(m, j) = \langle \hat{R}_n^m f, g_j \rangle_{(n)} - \langle \tilde{R}_n^m f, g_j \rangle$. Then, by definition,

$$\begin{aligned}
A_n(m, j) &= A_n(m-1, j) - \langle \tilde{R}_n^{m-1} f, g_{\hat{\mathcal{S}}_m} \rangle (\langle g_{\hat{\mathcal{S}}_m}, g_j \rangle_{(n)} - \langle g_{\hat{\mathcal{S}}_m}, g_j \rangle) \\
&\quad - \langle g_{\hat{\mathcal{S}}_m}, g_j \rangle_{(n)} (\langle \hat{R}_n^{m-1} f, g_{\hat{\mathcal{S}}_m} \rangle_{(n)} - \langle \tilde{R}_n^{m-1} f, g_{\hat{\mathcal{S}}_m} \rangle) \\
&= A_n(m-1, j) - I_{n,m}(j) - II_{n,m}(j).
\end{aligned} \quad (6.10)$$

From Lemma 1(i) we get

$$\sup_{1 \leq j \leq p_n} |I_{n,m}(j)| \leq \|\tilde{R}_n^{m-1} f\| \|g_{\hat{\mathcal{S}}_m}\| \zeta_n \leq \|f\| \zeta_n, \quad (6.11)$$

where we have used the norm-reducing property in (6.3) for $\tilde{R}_n^m f$.

For the second term we proceed recursively:

$$\begin{aligned}
\sup_{1 \leq j \leq p_n} |II_{n,m}(j)| &\leq \sup_{1 \leq j \leq p_n} |\langle g_{\hat{\mathcal{S}}_m}, g_j \rangle_{(n)}| \sup_{1 \leq j \leq p_n} |A_n(m-1, j)| \\
&\leq (1 + \zeta_n) \sup_{1 \leq j \leq p_n} |A_n(m-1, j)|.
\end{aligned} \quad (6.12)$$

For the last inequality, we have used again Lemma 1(i) and the Cauchy–Schwarz inequality $|\langle g_{\hat{\mathcal{S}}_m}, g_j \rangle| \leq \|g_{\hat{\mathcal{S}}_m}\| \|g_j\| = 1$.

Using the notation $B_n(m) = \sup_{1 \leq j \leq p_n} |A_n(m, j)|$, we get the following recursion from (6.10)–(6.12):

$$\begin{aligned}
B_n(0) &\leq \zeta_n, \\
B_n(m) &\leq B_n(m-1) + \zeta_n \|f\| + (1 + \zeta_n) B_n(m-1) \\
&\leq (5/2) B_n(m-1) + \zeta_n \|f\| \quad \text{on the set } A_n.
\end{aligned}$$

Therefore,

$$B_n(m) \leq (5/2)^m \zeta_n + \zeta_n \|f\| \sum_{j=0}^{m-1} (5/2)^j \leq (5/2)^m \zeta_n \left(1 + \|f\| \sum_{j=0}^{m-1} (5/2)^{j-m}\right)$$

$$\leq (5/2)^m \zeta_n \left(1 + \sup_{n \in \mathbb{N}} \sum_{j=1}^{p_n} |\beta_{j,n}| \sum_{k=1}^{\infty} (5/2)^{-k}\right),$$

which completes the proof by setting $C_* = 1 + \sup_{n \in \mathbb{N}} \sum_{j=1}^{p_n} |\beta_{j,n}| \sum_{k=1}^{\infty} (5/2)^{-k}$. □

*Analyzing $\tilde{R}_n^m f$.* We are now ready to establish a finite-sample analogue of (6.2) for $\tilde{R}_n^m f$. We have

$$\langle \hat{R}_n^m f, g_j \rangle_{(n)} = \langle \tilde{R}_n^m f, g_j \rangle + (\langle \hat{R}_n^m f, g_j \rangle_{(n)} - \langle \tilde{R}_n^m f, g_j \rangle).$$

Hence, by invoking Lemma 2 (and denoting by $A_n$ the set as there) we get

$$|\langle \hat{R}_n^m f, g_{\hat{\mathcal{S}}_m} \rangle_{(n)}| = \sup_{1 \leq j \leq p_n} |\langle \hat{R}_n^m f, g_j \rangle_{(n)}|$$

$$\geq \sup_{1 \leq j \leq p_n} |\langle \tilde{R}_n^m f, g_j \rangle| - (5/2)^m \zeta_n C_* \qquad \text{on the set } A_n.$$

Therefore, again by Lemma 2 for the first inequality to follow,

$$(6.13) \quad \begin{aligned} |\langle \tilde{R}_n^m f, g_{\hat{\mathcal{S}}_m} \rangle| &\geq |\langle \hat{R}_n^m f, g_{\hat{\mathcal{S}}_m} \rangle_{(n)}| - (5/2)^m \zeta_n C_* && \text{on the set } A_n \\ &\geq \sup_{1 \leq j \leq p_n} |\langle \tilde{R}_n^m f, g_j \rangle| - 2(5/2)^m \zeta_n C_* && \text{on the set } A_n. \end{aligned}$$

Consider the set $B_n = \{\omega; \sup_{1 \leq j \leq p_n} |\langle \tilde{R}_n^m f, g_j \rangle| > 4(5/2)^m \zeta_n C_*\}$. Then, by (6.13),

$$(6.14) \quad |\langle \tilde{R}_n^m f, g_{\hat{\mathcal{S}}_m} \rangle| \geq 0.5 \sup_{1 \leq j \leq p_n} |\langle \tilde{R}_n^m f, g_j \rangle| \qquad \text{on the set } A_n \cap B_n.$$

Formula (6.14) says that the selectors $\hat{\mathcal{S}}_m$ satisfy the condition (6.2) for $\tilde{R}_n^m f$ on the set $A_n \cap B_n$. We can now invoke Temlyakov's result in (6.5), since the condition (6.2) holds on the set $A_n \cap B_n$ [as established in (6.14)]. We have

$$(6.15) \quad \|\tilde{R}_n^m f\| \leq B(1 + m/4)^{-1/10} = o(1) \qquad \text{on the set } A_n \cap B_n$$

by choosing $m = m_n \to \infty$ ($n \to \infty$) (slow enough), where $B = \sup_{n \in \mathbb{N}} \sum_{j=1}^{p_n} |\beta_{j,n}| < \infty$; see (6.4) and assumption (A2).

For $\omega \in B_n^C = \{\omega; \sup_{1 \leq j \leq p_n} |\langle \tilde{R}_n^m f, g_j \rangle| \leq 4(5/2)^m \zeta_n C_*\}$, by using formula (5.2) from [21] with $b_m$ as defined there (i.e., $b_m = b_{m-1} + |\langle \tilde{R}_n^{m-1} f, g_{\hat{\mathcal{S}}_m} \rangle|$, $b_0 = 1$),

$$(6.16) \quad \begin{aligned} \|\tilde{R}_n^m f\|^2 &\leq \sup_{1 \leq j \leq p_n} |\langle \tilde{R}_n^m f, g_j \rangle| b_m \\ &\leq \sup_{1 \leq j \leq p_n} |\langle \tilde{R}_n^m f, g_j \rangle|(1 + m\|f\|) \\ &\leq 4(5/2)^m \zeta_n C_*(1 + m\|f\|) \qquad \text{on the set } B_n^C. \end{aligned}$$

For bounding the number $b_m$, we have used the norm-reducing property in (6.3) applied to $\tilde{R}_n^m f$. Therefore, using (6.15), (6.16) and $\zeta_n = O_P(n^{-\xi/2})$ from Lemma 1, we have for $m = m_n \to \infty$ ($n \to \infty$) slow enough [e.g., $m_n = o(\log(n))$],

$$(6.17) \quad \begin{aligned} \|\tilde{R}_n^m f\| &\leq B(1 + m_n/4)^{-1/10} \\ &\quad + 4(5/2)^{m_n} \zeta_n C_*(1 + m\|f\|) \qquad \text{on the set } (A_n \cap B_n) \cup B_n^C \\ &= o_P(1), \end{aligned}$$

since $\mathbb{P}[(A_n \cap B_n) \cup B_n^C] \geq \mathbb{P}[A_n] \to 1$, $n \to \infty$, due to Lemma 1.

*Analyzing* $\hat{R}_n^m f$. By definition and using the triangle inequality,

(6.18) $\qquad \|\hat{F}_n^m - f\| = \|\hat{R}_n^m f\| \le \|\tilde{R}_n^m f\| + \|\hat{R}_n^m f - \tilde{R}_n^m f\|.$

A recursive analysis can be developed for the second term on the right-hand side:

$$A_n(m) = \|\hat{R}_n^m f - \tilde{R}_n^m f\|.$$

By definition,

$$A_n(m) = \|\hat{R}_n^{m-1} f - \tilde{R}_n^{m-1} f - (\langle \hat{R}_n^{m-1} f, g_{\hat{S}_m}\rangle_{(n)} - \langle \tilde{R}_n^{m-1} f, g_{\hat{S}_m}\rangle) g_{\hat{S}_m}\|$$

$$\le A_n(m-1) + |\langle \hat{R}_n^{m-1} f, g_{\hat{S}_m}\rangle_{(n)} - \langle \tilde{R}_n^{m-1} f, g_{\hat{S}_m}\rangle| \|g_{\hat{S}_m}\|$$

$$\le A_n(m-1) + (5/2)^{m-1} \zeta_n C_* \qquad \text{on the set } A_n,$$

where the last inequality follows from Lemma 2. Therefore, for some constant $C > 0$,

(6.19) $\qquad \|\hat{R}_n^m f - \tilde{R}_n^m f\| \le 3^m \zeta_n C = o_P(1)$

by choosing $m = m_n \to \infty$ sufficiently slowly such that $3^{m_n} \zeta_n = o_P(1)$.

By (6.17)–(6.19) we get [e.g., by using the choice $m_n \to \infty, m_n = o(\log(n))$]

$$\mathbb{E}_X |\hat{F}_n^{m_n}(X) - f(X)|^2 = \|\hat{R}_n^{m_n} f\|^2 = o_P(1),$$

which completes the proof of Theorem 1.

6.3. *Arbitrary step-size* $\nu$. For arbitrary, fixed step-size $0 < \nu \le 1$ in (2.2), we need to make a few modifications to the proof.

Temlyakov's result in (6.5) becomes

$$\|R^m f\| \le B(1 + \nu(2-\nu) m b^2)^{-b/(2(2+b))}, \qquad 0 < b \le 1 \text{ as in } (6.2).$$

PROOF. The claim follows as in [21]. Using his notation, we use $a_m = \|R^m f\|^2$, $y_m = |\langle R^{m-1} f, g_{S_m}\rangle|$, $b_m = b_{m-1} + \nu y_m, b_0 = 1$ and $t_m \equiv b$ from (6.2). □

We can then use exactly the same reasoning as in Section 6.2. At some obvious places, a factor $\nu$ occurs in addition, and it can be trivially bounded by 1. The only slightly nontrivial reasoning occurs in (6.16); but using $b_m$ as defined above (applied now to $\tilde{R}_n^m f$ instead of $R^m f$) yields the bound

$$\|\tilde{R}_n^m f\|^2 \le \sup_{1 \le j \le p_n} |\langle \tilde{R}_n^m f, g_j\rangle| b_m,$$

which then allows us to proceed as in Section 6.2.

6.4. *Binary classification.* The first assertion of Corollary 1 follows exactly as in the proof of Theorem 1 by using the representation in (3.2). There is no crucial place where we make use of homoscedastic errors $\varepsilon_i$: the uniform laws of large numbers from Lemma 1 look formally a bit different (e.g., for (ii) we need to subtract a term $\mathbb{E}[g_j(X_i)\varepsilon_i] = \mathbb{E}[g_j(X_i)(Y_i - f(X_i))]$), but the i.i.d. structure of the pairs $(X_i, Y_i)$ suffices to get through. The moment assumption for $\varepsilon_i = Y_i - f(X_i)$ trivially holds since $|Y_i| \leq 1$ and $\sup_x |f(x)| \leq 1$.

For the second assertion, it is well known that 2 times the $L_1$-norm bounds from above the difference between the generalization error of a plug-in classifier (expected 0–1 loss error for classifying a new observation) and the Bayes risk ([8], Theorem 2.3). Furthermore, the $L_1$-norm is upper bounded by the $L_2$-norm.

## APPENDIX

**The model (4.4).** The model (4.4) is as follows. Define $a_j = j^{0.51}$. Let the parameter $\kappa$ be the solution of the equation $\sigma_\varepsilon^2 n^{-1} \sum_{j=1}^\infty a_j \lambda_j = \kappa$, where we denote $\lambda_j = (1 - \kappa a_j)_+$. For $n = 100$, the solution is $\kappa = 0.199$. Determine the predictor dimension $p = \max_j \{a_j \leq \kappa^{-1}\} = 23$. The variances are $\sigma_j^2 = \lambda_j (n\kappa a_j)^{-1}$, $j = 1, \ldots, 23$, $n = 100$. It can be shown that such regression coefficients belong with high probability to $\{(\beta_{j,n})_j; \sum_{j=1}^{p_n} a_j^2 \beta_{j,n}^2 \leq 1\}$ (note that $p = p_n$ depends on $n$ via the parameter $\kappa = \kappa_n$).

**Acknowledgments.** I thank two referees, an Associate Editor and the Editor M. Eaton for constructive comments.

SEMINAR FÜR STATISTIK
ETH ZÜRICH
CH-8092 ZÜRICH
SWITZERLAND
E-MAIL: buhlmann@stat.math.ethz.ch